\documentclass[twoside, 11pt]{article}

\usepackage[latin1]{inputenc}
\usepackage{xspace}
\usepackage{amsmath,amssymb,amsfonts,epic}
\usepackage[dvips]{graphicx}
\usepackage{cite}
\usepackage{graphics,epsfig}
\usepackage{fancybox}
\usepackage{amsbsy}
\usepackage{amscd}
\usepackage{amsgen}
\usepackage{amsopn}
\usepackage{amstext}
\usepackage{ifthen}
\usepackage{amsxtra}
\usepackage{epsfig}
\usepackage{color}
\usepackage{fancyhdr}
\usepackage{fullpage}
\usepackage{float}
\usepackage[nottoc]{tocbibind}

\definecolor{hellgrau}{gray}{0.4}

\newcommand{\e}{\varepsilon}
\newcommand{\Ep}{\mathcal{E}}
\newcommand{\N}{\mathbb{N}}
\newcommand{\Z}{\mathbb{Z}}

\newcommand{\E}{\mathbb{E}}

\newcommand{\1}{1\! \mathrm{l}}

\newcommand{\y}{\mathrm{y}} 

\newtheorem{theorem}{Theorem}

\newtheorem{corollary}{Corollary}

\newtheorem{exist.cond}{Existence conditions}
\newtheorem{remark}{Remark}
\newtheorem{lemma}{Lemma}
\newenvironment{proof}[1][Proof] {{\bf{#1}}}{\hfill $\square$}

\begin{document}

\title{Aggregation of autoregressive processes and long memory}

\author{By   Didier Dacunha-Castelle$^1$ and Lisandro J. Ferm\'{\i}n$^{1,2}$}
\date{}
\maketitle

\begin{center}
\begin{minipage}{0.6\textwidth}
\begin{center}
{\tiny
$^1$ UMR CNRS 8628. \'Equipe de Probabilités, Statistique et Modélisation, Bât. 425,\\
Université Paris-Sud, 91405 Orsay Cedex, France.\\
\vspace{0.25cm}
$^2$ UMR CNRS 8145 Laboratoire de Mathématiques Appliquées à Paris 5 (MAP5),\\
Université Paris Descartes, 45 rue des Saints Pères 75270 Paris Cedex 06, France.\\

e-mail: didier.dacunha-castelle@math.u-psud.fr, fermin.lisandro@math.u-psud.fr
}
\end{center}
\end{minipage}
\end{center}
\vspace{0.2cm}

\abstract{
\vspace{0.2cm}
\begin{center}
\begin{minipage}{0.8\textwidth}
\hspace{0.4cm} We study the aggregation of AR processes and generalized Ornstein-Uhlenbeck (OU) processes. Mixture of spectral densities with random poles are the main tool. In this context, we apply our results for the aggregation of doubly stochastic interactives processes, see \cite{Dacunha&Fermin.L.Notes}. Thus, we study the relationship between aggregation of autoregressive processes and long memory considering complex interaction structures. We precise a very interesting qualitative phenomena: how the long memory creation depends on the poles concentration near to the boundary of stability (measured in the Prokhorov sense).  Our results extends the results given by Oppenheim and Viano, \cite{Oppenheim&Viano2004}, and highlight the importance of the angular dispersion measure of poles in the appearance of the long memory.

\end{minipage}
\end{center}
\vspace{0.2cm}
}
\thispagestyle{empty}


\section{Introduction}

\hspace{0.4cm} Long memory (LM) processes, are used in many fields such as economics, finance,  hydrology  or  communication  networks.  Some of  these  LM processes can be  seen as an aggregation of elementary short memory (SM) processes.

The aggregation of stochastic processes was introduce by Granger in 1980, \cite{Granger}. It is a sommation procedure of identically distributed elementary processes $Z^i=\{Z^i_t \, : t\in \Z\}$ over the index $i$. Granger shows  that  by aggregating  random parameter $AR(1)$ processes, one can obtain LM processes with spectral density equivalent   to  $|\lambda|^{-d}$ when $\lambda\rightarrow 0$
for some  $d$, $0<d<1$. He  considers $AR(1)$ processes with independent random Beta distributed parameters and gives conditions on  the Beta distribution in order to obtain long memory. That has opened a new way of obtaining long range dependence time series.

As indicated by Beran, \cite{Beran}, this is  an interesting idea from two different points of view: it allows a physical explanation
of the long dependence phenomenon in several fields  and it gives an easy and fast simulation method.

We develop the procedure of aggregation considering doubly stochastic $AR(p)$ (or generalized Ornstein-Uhlenbeck process of order $p$) elementary processes $Z^i= \{Z_t(\y^i, \varepsilon^i): \, t \in T\}$, where $Y=\{\y^i: \, i \in \mathbb{N}\}$ is a sequence of random
variables with distribution $\mu$ on $\mathbb{R}^p$ and $\Ep= \{\varepsilon^i: i \in \N\} $ is the sequence of innovations $\e^i= \{\varepsilon^i_t : \, t \in T\}$. Finally, $T= \mathbb{Z}$ in the case of discrete time processes and $T \subseteq \mathbb{R}$ in the
continuous time case. We assume that $Y$ and $\Ep$ satisfy the following assumption.

\quad\\
\textbf{Assumption A1:}
\begin{enumerate}
  \item  \textit{$\Ep$ is an array of strong white noises, i.e. for each $i$, $\e^i=\{\e^i_t\}$
  is an i.i.d sequence such that $\E[|\varepsilon^i_{t}|^2]=1$.}

  \item \textit{$Y$ is an i.i.d. sequence with distribution $\nu= \mu^{\bigotimes \N}$.}

  \item  \textit{$Y$ is independent of $\Ep$.}
\end{enumerate}

Then, we define the sequence of partial aggregations $X^N(Y)=\{X^N_t(Y):\, t \in T \}$ of
elementary processes $\{Z^i\}$, by
\begin{equation}\label{ec.7.1.0}
X^N_t(Y)=\frac{1}{B_N}\sum_{i=1}^N Z_t(\y^i, \varepsilon^i),
\end{equation}
where $\{B_N\}$ is a normalization sequence.

In general, we consider three types of innovations:
\begin{enumerate}

\item \textbf{Common innovation:} $\e^i= \e$ for all $i \in \N$.

\item \textbf{Independent innovations:} $\e^i$ independent of $\e^j$ for $i \neq
j$.

\item \textbf{Stationary interactive innovations:}  $\mathbb{E}[\varepsilon^i_t   \,
\varepsilon^j_t   ]=\chi(i-j)$   and   $\mathbb{E}[\varepsilon^i_t \, \varepsilon^j_s ]=0$ for
$t\neq s$, where $\chi$ is an interaction correlation.
\end{enumerate}

The aggregation procedure of doubly stochastic linear processes can be also develop using mixtures of spectral densities and mixture of transfer functions as main tools, see \cite{Dacunha&Fermin.L.Notes}. Let $g(\lambda,\y)$, $ \y \in \mathbb{R}^p$, be a family of spectral densities and $\mu$ be a probability on $\mathbb{R}^p$. We will denote by $h(\cdot, \y )$ a particular square root of $g(\cdot, \y)$, some time for simplicity we consider the real one. In other case, there exists a root linked with a regular representation. This is always possible for autoregressive process of order $p$.

We consider, for $\y^i$ fixed, that $g(\lambda,\y^i)$  is  the spectral density  of the elementary process $Z^i=Z(\y^i, \e^i)$. The mixture of the spectral densities $g(\lambda,\y)$, is defined by
\begin{equation} \label{ec.7.1.1}
F(\lambda)=\int_{\mathbb{R}^p} g(\lambda,\y) \mu(d \y).
\end{equation}
$F(\lambda)$ is a well defined spectral density if and only if
\begin{equation}\label{ec.7.1.2}
\int F(\lambda)d\lambda < \infty.
\end{equation}

The mixture of transfer function given by
\begin{equation}\label{ec.7.1.3}
H(\lambda)=\int_{\mathbb{R}^s}h(\lambda,\y) \mu(d\y) \; ,
\end{equation}
is well defined and will be  called a transfer function if $|H|^2$   is   a   spectral density. This is a consequence of condition~\eqref{ec.7.1.2}  and  Jensen's  inequality.

Without loss of generality for the purpose of this chapter, we assume that $\Ep=\{\e^i\}$ is a
sequence of Gaussian white noise, except for some remarks.

Under some general conditions for the interaction $\chi$ we show in \cite{Dacunha&Fermin.L.Notes} the $\nu-a.s.$ weak convergence of $X^N(Y)$ to a process $X$ called the aggregation process.

In the independent innovations case, we take $B_N = \sqrt{N}$, then one can show that the aggregation process exists $\nu-a.s.$ if and only if condition~\eqref{ec.7.1.2} hold. In this case the aggregation process $X$ is Gaussian with spectral density $F$. Nevertheless, for the common innovation case, taking $B_N = N$ we can prove that condition~\eqref{ec.7.1.2} implies the existence $\nu-a.s.$ of $X$, in this case $X$ is a Gaussian process with spectral density $|H|^2$.

For the case of interactive innovations we have show, in \cite{Dacunha&Fermin.L.Notes}, that in general condition~\eqref{ec.7.1.2} is a necessary and sufficient condition for the existence $\nu-a.s.$ of aggregation process $X$.  Moreover, the limit is always a convex combination of the two extreme cases: independent innovations ($F(\lambda)$, $B_N = \sqrt{N}$) and common innovations ($H^2(\lambda)$, $B_N= N$). In general, the limit is reached for a normalization $B_N$ which depends on the behavior of interaction $\chi$. Ferm\'{\i}n, \cite{Fermin}, generalize this results for the case of non-gaussian innovations, considering a sequence $\Ep=\{\e^i\}$ of weakly dependent innovations.

Thee long memory of a stationary process is de defined by the nonsummability of the corresponding covariance sequence $\gamma=\{\gamma(k): \, k \in \Z\}$; i.e. when the norm $\|\gamma\|_{1}= \sum_k |\gamma(k)|<\infty$ we say that the process is SM and if $\|\gamma\|_{1}= \infty$ then we say that the process is LM. The long memory sequences are generally associated to the singularities of the spectral density.

The non-uniformity  of $\|\gamma(\cdot, \y)\|_{1}$, for different values of $\y$, can generate long memory by aggregation. In the spectral density context, this phenomenon is related with the concentration of the mixture measure $\mu$ near the boundary of the existence domain for elementary processes, defined  by the set of parameters $\y$ such that $g(\lambda, \y)$ is a well defined spectral density; i.e. the set $\{\y: \, \int g(\lambda,\y) d\lambda <  \infty\}$.  For instance, in the case of $AR(1)$ processes aggregations, we  have an only way of creating long  memory; that is, when the measure $\mu$ of the autoregressive parameter is concentrated near
$\{-1,1\}$, which is the boundary of the domain of parameter values leading to stationary $AR(1)$.

Several authors,  \cite{Goncalvez1988, Lippi&Zaffaroni, Linden, Terence1, Celov}, study the long memory on the aggregation of $AR(1)$ and $AR(2)$ processes, but considering only the case in which the spectral density  has an only singularity at zero, which induces long memory but without seasonal effects, for the continuous time case see \cite{Barndorff, Igloi&Terdik, Nikolai}. In \cite{Oppenheim&Viano2004}, the authors show how to obtain seasonally LM models from the aggregation of autoregressive processes in the discrete time case as well as in the continuous time case.

Our purpose is to clarify in a more general context, the conditions linked to the existence and the long memory property  of the aggregation of $AR$ and $OU$ processes. The results given in this paper  are an extension of the results shown in \cite{Oppenheim&Viano2004}. Our generalization consists in considering poles with multiplicity and with diffuse angular distribution.

We shall exhibit a new qualitative behavior: for $p\geq 2$ the long memory is not guaranteed even if the radial distribution of poles diverges in $1$; the angular distribution of poles is also determinant and can counterbalance the radial distribution. For instance, for $AR(2)$ elementary processes, if $(\rho e^{i\theta})^{-1}$ is a random pole and if $\theta$ has a diffuse measure then $\rho$ has to be very concentrated near $\{-1,1\}$ whereas if $\theta$ has a concentrated measure near to a given frequency then conditions on $\rho$ can be relaxed.

Results are presented as follows. In Section~\ref{sect.7.2}, we show that for suitable distributions of the $AR(p)$ parameters the aggregation process exists  and we give the expression of its spectral density. Furthermore,  we give necessary  and sufficient conditions for this spectral density to have predetermined singularities which imply the long memory, we study in detail the cases of $AR(1)$ and $AR(2)$ processes and we present a general result for the case $AR(p)$. We illustrate how the long memory can "disappear" when the angular distribution of of poles is diffuse. In Section~\ref{sect.7.4} we extend the results given in Section~\ref{sect.7.2} to the case of the aggregation of $OU(p)$ processes.


\section{Aggregation of $AR(p)$ processes and long memory \label{sect.7.2}}

\hspace{0.4cm} Let  $D$ be the open unit disc. For  $ \y = (\y_1,..., \y_p)  \in D^p$, let $Z_t(\y, \e)$ denote  the autoregressive process  with innovation $\varepsilon$ satisfying
\begin{equation}\label{ec.7.2.0}
Z_t = - \sum_{k=1}^p a_k(\y) Z_{t-k} + \sigma \e_t.
\end{equation}
We assume that the corresponding characteristic polynomial is factorized of the following way
$$A(s,\y)=1+ \sum_{k=1}^p a_k(\y) s^k = \prod_{k=1}^p(1-\y_k s).$$
So  $\{1/\y_k: \, k=1,...,p\}$ is the set of  random poles of the $AR(p)$ process. The process $Z_t(\y, \e)$ has the $MA(\infty)$ expansion
\begin{equation}\label{ec.7.2.1}
Z_t(\y)= \sigma \sum_{k=0}^{\infty}c_k(\y)\varepsilon_{t-k},
\end{equation}
where $ A^{-1}(s,\y)=\sum_{k=0}^{\infty}c_k(\y)s^k$.

Suppose that $\y$ is a  random  vector  independent of the  innovation $\varepsilon$ and whose
distribution $\mu$ has support $D^p$. From  the independence assumption  and  because
$\mu(|\y_k|<1)=1$ for every $k$, the  series in \eqref{ec.7.2.1} converge almost  surely  for
$\mu$-almost all $\y$. The induced process $Z_t(\y, \e)$ is, for almost all $\y$, a stationary $AR(p)$
process with spectral density $g(\lambda,\y)=\sigma^2 |A(e^{i\lambda},\y)|^{-2}$. We take the
transfer function $h(\lambda,\y)=\sigma A(e^{i\lambda}, \y)^{-1}= \sigma
\sum_{k=0}^{\infty}c_k(\y)e^{i\lambda k}$.

Since $F(\lambda)= \E[g(\lambda,\y)]$, $H(\lambda)= \E[h(\lambda,\y)] $ and
$$\E^{\y}[|Z_t(\y)|^2]=\int_{-\pi}^{\pi}g(\lambda,\y)d\lambda,$$
where $\E^{\y}[\cdot]$ denotes the conditional expectation given $\y$, then we have that $Z_t(\y, \e)$
exists in $L^2$ if and only if condition \eqref{ec.7.1.2} holds.

We take $\{Z^i\}$ as a sequence of random parameters $AR(p)$ processes, with $Y=\{\y^i\}$ a
sequence of random vectors in $D^p$. Without loss of generality for the purpose of this paper, we
assume that $\Ep=\{\e^i\}$ is a sequence of Gaussian white noise. We suppose that $Y$ and $\Ep$
satisfy the Assumption~A1 and so we have the convergence results of Theorem 1 given in \cite{Dacunha&Fermin.L.Notes}.

We denote by $\mathbb{AR}(p)$ the class of $AR(p)$ processes and by $\mathbb{M}(\mathbb{AR}(p))$ the class of processes that can be obtained by aggregation of elementary processes in $\mathbb{AR}(p)$.

In this section, our aim is to give a way of taking the measure $\mu$ that allows us to obtain the long memory property for $X$.

From now, we will consider that $A(s,\y)$ is a polynomial of degree $p$ with  $n$ real roots  and $2(q-n)$  complex pairwise  conjugate roots; i.e.,
$$A(s,\y)=\prod_{k=1}^n(1 \pm \rho_k s)^{m_k}\prod_{k=n+1}^{q}\left[(1-\rho_k e^{i\theta_k}s)(1-\rho_k e^{-i\theta_k}s)\right]^{m_k},$$
where $n\leq  q$, $m_k$ is  the multiplicity of the  root $\y_k^{-1}$, where $\y_k=\rho_k e^{i\theta_k}$ with $\theta_k \in [-\pi,\pi)$, $\rho_k \in (0,1)$ and $p=\sum_{k=1}^{n}m_k+ 2\sum_{j=n+1}^{q}m_j$. Furthermore, we will consider  the following assumption.

\quad\\
\textbf{Assumption B1:} \textit{Let $\rho_1,...,\rho_q$ be  independent random variables, such that
$\rho_k$ has for distribution
\begin{equation}\label{ec.7.3.1}
dR_k(\rho)=|1-\rho|^{d_k} \varphi_k(\rho)d\rho,
\end{equation}
where $\varphi_k$  is a bounded  positive function with support $[0,1]$, continuous  in   $\rho=1$
with  $\varphi_k(1)>  0$. Let $\theta_{1},..., \theta_q$ be independent random variables and
independent of $\rho_1,...,\rho_q$, such that $\theta_k$ has for distribution
\begin{equation}\label{ec.7.3.2}
dQ_k(\theta)=\frac{\psi_k(\theta)d\theta}{|\theta-\theta_k^0|^{\beta_k}},
\end{equation}
where  $\beta_k  \leq  1$,  $\psi_k$  is  a  bounded  positive function with support  $(-\pi,\pi]$,
continuous  in $\theta=\theta_k^0$ and $\psi_k(\theta_k^0)>0$.}

\begin{remark}\label{rm.7.3.1}
Note that when $\beta_k \rightarrow 1$ and $\psi_k(\theta)  \rightarrow 0$  for $\theta  \neq
\theta_k^0$, then the measure $Q_k$  converges   to  Dirac's  delta;  i.e., by convention  when
$\beta_k=1$               we               will consider $Q_k(\theta)=\delta(\theta-\theta_k^0)$.
Thus, we consider $\beta_k = 1$  and $\theta_k^0  \in \{0,\pi \}$,  for $1\leq k\leq  n$, without
more details.
\end{remark}

Then, if condition \eqref{ec.7.1.2} is fulfilled we have that
\begin{equation}\label{ec.7.3.3}
F(\lambda)=\prod_{k=1}^n\int_0^1\frac{dR_k(\rho)}{|1 - \rho e^{i(\lambda- \theta_k^0 )}|^{2m_k}}
\prod_{k=n+1}^q \int_0^1\int_{-\pi}^{\pi}\frac{dR_k(\rho)dQ_k(\theta)}{[|1-\rho
e^{i(\lambda+\theta)}||1-\rho e^{i(\lambda-\theta)}|]^{2m_k}}.
\end{equation}
\begin{equation}\label{ec.7.3.4}
H(\lambda)=\prod_{k=1}^n\int_0^1\frac{dR_k(\rho)}{(1- \rho e^{i(\lambda-\theta_k^0)})^{m_k}}
\prod_{k=n+1}^q \int_0^1\int_{-\pi}^{\pi}\frac{dR_k(\rho)dQ_k(\theta)}{[(1-\rho
e^{i(\lambda+\theta)})(1-\rho e^{i(\lambda-\theta)})]^{m_k}}.
\end{equation}

\begin{remark}\label{rm.7.3.2}
We study in detail the cases of independent innovations and common innovation. Since, in the case of interactive innovations the spectral density of the aggregated process is always a positive convex combination of the form $a F + b|H|^2$ then the results in this last case can be deduced from the two previous cases.
\end{remark}

In the following we will consider that  $A_{\lambda} \sim B_{\lambda}$ near $\lambda  = \lambda_0$, if  $\lim_{\lambda \rightarrow \lambda_0} A_{\lambda}/ B_{\lambda}$ is a non-null constant. Under
Assumption~B1 we have the following lemmas.

\begin{lemma}\label{lemma01} \quad
\begin{enumerate}

\item If $-1<d_k< n_k-1$, then near $\lambda=0$
\begin{eqnarray*}
\int_0^1\frac{|1-\rho|^{d_k}\varphi_k(\rho)d\rho}{|1-\rho e^{i\lambda}|^{n_k}}  \sim
\frac{\varphi_k(1)}{|\lambda|^{n_k-1-d_k}}\int_0^{\infty}\frac{u^{d_k}du}{(1+u^2)^{n_k/2}}.
\end{eqnarray*}

\item If $-1<d_k<n_k-1$, then near $\lambda=\pi$
\begin{eqnarray*}
\int_0^1\frac{|1-\rho|^{d_k}\varphi_k(\rho)d\rho}{|1+\rho e^{i\lambda}|^{n_k}}
 \sim\frac{\varphi_k(1)}{|\lambda-\pi|^{n_k-1-d_k}}\int_0^{\infty}\frac{u^{d_k}du}{(1+u^2)^{n_k/2}}.
\end{eqnarray*}

\item If $-1<d_k<n_k-1$ and  $\theta_k^0 \notin \{0,\pi\}$, then near
$\lambda =\pm \theta_k^0$
\begin{eqnarray*}
\int_0^1\!\!\frac{|1-\rho|^{d_k}\varphi_k(\rho)d\rho}{[|1-\rho e^{i(\lambda+\theta_k^0)}||1-\rho
e^{i(\lambda-\theta_k^0)}|]^{n_k}}\sim \frac{\varphi_k(1)[2sin(\theta_k^0)]^{-n_k}}
{|\lambda\mp\theta_k^0|^{n_k-1-d_k}}\!\!\int_0^{\infty}\!\!\!\frac{u^{d_k}du}{(1+u^2)^{n_k/2}}.
\end{eqnarray*}
\end{enumerate}
\end{lemma}

\begin{lemma}\label{lemma02}Let $\alpha_k<1$ and
$$f_k(\lambda)=\int_{-\pi}^{\pi} \int_0^1 \frac{|1-\rho|^{d_k}\varphi_k(\rho)\psi_k(\theta)
|\theta-\theta_k^0|^{-\alpha_k} d\rho d\theta} {[|1-\rho e^{i(\lambda+\theta)}||1-\rho e^{i(\lambda-\theta)}|]^{n_k}}.$$
\begin{enumerate}
\item If $ n_k-1<d_k<2n_k-2+\alpha_k$, $\theta_k^0=0$, then near $\lambda = 0$
\begin{eqnarray*}
f_k(\lambda) \sim \frac{\varphi_k(1)\psi_k(0)}{|\lambda|^{2n_k-2-d_k+\alpha_k}}
\int_{-\infty}^{\infty}\!\int_0^{\infty}\!\! \frac{u^{d_k}|\theta|^{-\alpha_k}du
d\theta}{[((\theta-1)^2+u^2)((\theta+1)^2+u^2)]^{n_k/2}}.
\end{eqnarray*}

\item If $ n_k-1<d_k<2n_k-2+\alpha_k$, $\theta_k^0=\pi$, then near
$\lambda = \pi$
\begin{eqnarray*}
f_k(\lambda) \sim \frac{\varphi_k(1)\psi_k(\pi)}{|\lambda-\pi|^{2n_k-2-d_k+\alpha_k}} \!
\!\int_{-\infty}^{\infty} \!  \int_0^{\infty}\!\! \frac{u^{d_k}|\theta|^{-\alpha_k}du
d\theta}{[((\theta-1)^2+u^2)((\theta+1)^2+u^2)]^{n_k/2}}.
\end{eqnarray*}

\item If $n_k-2<d_k<n_k-2+\alpha_k$ and $\theta_k^0\notin \{0,\pi\}$, then near $\lambda = \pm \theta_k^0$
\begin{eqnarray*}
f_k(\lambda) \sim \frac{\varphi_k(1)\psi_k(\theta_k^0)[2sin(\theta_k^0)]^{-n_k}}{|\lambda \mp
\theta_k^0|^{n_k-2-d_k+\alpha_k}}
\int_{-\infty}^{\infty}\!\!\frac{|\theta|^{-\alpha_k}d\theta}{|\theta \mp 1|^{n_k-1-d_k}}
\int_0^{\infty}\!\!  \frac{u^{d_k} du}{(1+u^2)^{n_k/2}}.
\end{eqnarray*}

\end{enumerate}
\end{lemma}

We do not give the proof of these lemmas, since they are similar to those given for Lemma~\ref{lemma1} and Lemma~\ref{lemma2} for the continuous time case, we referred to Section~\ref{sect.7.4}.

When the measures $dR_k$ are concentrated near of the boundary
$$\delta D^p= \left\{\y: \sup_{1\leq k \leq p}|\y_k|=1\right\}$$
of $D^p$, then the first $n$ terms  in \eqref{ec.7.3.3}, or in \eqref{ec.7.3.4}, can only produce a singularity in $F$, or respectively  in  $H$, at  the frequencies $0$ or $\pi$, while singularities at other frequencies can be provided by the last terms. So, for the study of the long memory, the behavior of the measures $dR_k$  near the boundary $\delta D $ is essential.

If the mixture probabilities $dQ_k(\theta)$ are regular or very diffuse, for instance the Lebesgue
measure on some finite interval, and if their supports do not intersect $\{0,\pi\}$, then the  long
memory induced by the  $\rho$ concentration  near $1$ can "disappear"; i.e., it  is not enough to
have the  mixture probabilities $dR_k$ concentrated near $\delta D^p $ to  reach LM by aggregation
of the random parameters $AR(p)$ processes. In  fact when the aggregation process exists and $dR_k$
are concentrated near $\delta D^p $ then is sufficient  that the probabilities $dQ_k$ are close to
probabilities with support of Lebesgue measure 0. This is a new result.

In the following we characterize the probability measures $dR$ and $dQ$ in order to make condition
\eqref{ec.7.1.2} hold and to obtain the long memory property of aggregation process $X$. First, we
study in detail the case $p=1$ and $p=2$, and then we give an example in the case $p=2$ where the
long memory "disappears" by randomness of the parameter $\theta$. Finally, we present the general
result in the case of $AR(p)$ processes.


\subsection{Case of $AR(1)$ processes \label{sect.7.3.1}}

\hspace{0.4cm} In this section we study the aggregation of random parameter $AR(1)$ processes considering dependence between individual innovations in order to show the influence of interactive innovations on the construction of LM processes.

From convergence results given in \cite{Dacunha&Fermin.L.Notes}, we have that a necessary and sufficient condition to obtain the existence of aggregation process is that the interaction correlation $\chi$ has a limit in the Cesaro sense. In this case the limit $s$ is such that $-\frac{1}{2}\leq s \leq \infty$. Thus we consider the following two types of interactions.

\begin{itemize}
\item Weak interaction: when $-\frac{1}{2}\leq s <\infty$. For instance, short interaction such that $\sum_j |\chi(j)|<\infty$, or large range moderate oscillation when $\sum_j \chi(j)<\infty$ and $\sum_j |\chi(j)|=\infty$.

\item Strong interaction: when $s =\infty$. For instance, when $\sum_j \chi(j)<\infty$.
\end{itemize}

We have that the spectral density $F$ and the transfer function $H$ are given by
\begin{eqnarray*}
F(\lambda) & = & \int_{0}^{1} \frac{\sigma^2}{|1-\rho e^{i(\lambda-\theta_0)}|^2}dR(\rho) \; .\\
H(\lambda) & = & \int_{0}^{1} \frac{\sigma}{(1-\rho e^{i(\lambda-\theta_0)})} dR(\rho)\; .
\end{eqnarray*}

When $dR$ is concentrated near enough the boundary $\delta D=\{1\}$ of $D$, we can produce a singularity on $F$ and on $H$ at the frequencies $\theta_0 \in \{0, \pi\}$. By taking, $dR(\rho)$ as in \eqref{ec.7.3.1} and applying Lemma~\ref{lemma01} we can verify that
\begin{enumerate}
\item [i.] If $-1<d<1$, then near $\lambda = \theta_0$ $F(\lambda) \cong \frac{1}{|\lambda - \theta_0|^{1-d}}.$
\item [ii.]If $-1<d<0$, then near $\lambda = \theta_0$ $|H(\lambda)|^2 \cong \frac{1}{|\lambda - \theta_0|^{-2d}}.$
\end{enumerate}

Furthermore, we have that $\int F(\lambda) d\lambda < \infty$ if and only if $d>0$ and $\int |H|^2(\lambda) d\lambda < \infty$ if and only if $d> -\frac{1}{2}$. Then we obtain the following theorem.

\begin{theorem}\label{theo.7.3.1}[ Aggregation of $AR(1)$ processes and long memory.]
If we consider the aggregation of $AR(1)$ processes with random parameter $\y$ satisfying Assumption~B1 , then we have
\begin{enumerate}
\item \textbf{Independent innovation case:} the aggregation $X$ exists if and only if $d>0$ and $X$ is a long memory process if and only if $d<1$.

\item \textbf{Common innovation case:} the aggregation $X$ exists if and only if $d>-\frac{1}{2}$ and it is a long memory process if and only if $d<0$.

\item \textbf{Interactive innovation case:}
\begin{enumerate}
\item [3.1.] Weak interaction: the aggregation $X$ exists if and only if $d>0$ and it is a long memory process if and only if $d<1$. In this case $|H|^2$ does not produce long memory.
\item [3.2.] Strong interaction: we obtain the same result that for common innovation.
\end{enumerate}
\end{enumerate}
\end{theorem}
\vspace{0.3cm}

From the above result follow two qualitative ways of obtaining $\alpha$-LM processes, for $0<\alpha<1$; i.e. LM processes with spectral density $G(\lambda)$ such that $G(\lambda)\sim \frac{1}{|\lambda- \theta_0|^{\alpha}}$ near $\lambda=\theta_0$:
\begin{enumerate}
\item [1.] If $0<d<1$ then considering weak interaction between innovations, we can obtain by aggregation $\alpha$-LM processes with $0\!<\!\alpha \!< \!1$ from $F$ contribution. In this case $|H|^2$ does not produce LM.

\item [2.] If $-\frac{1}{2}< d <0$ then considering strong long interaction between innovations, we can also obtain by aggregation $\alpha$-LM processes with $0<\alpha<1$, from $|H|^2$ contribution but for a much stronger concentration of the mixture measure near $\delta D$.
\end{enumerate}
\vspace{0.3cm}


\subsection{Case of $AR(2)$ processes \label{sect.7.3.2}}

\hspace{0.4cm} In this section we give a complete analysis of aggregation of $AR(2)$ under Assumption~B1. We consider two type of $AR(2)$ processes. The first type are $AR(2)$ processes with different real poles $\y_1=\rho_1$, $\y_2=\rho_2$. The second type are $AR(2)$ processes with complex conjugated random poles $\y_1=\rho e^{i\theta}$, $\y_2=\rho e^{-i\theta}$ (for which we obtain the particular case of doubly real poles when $\theta \in \{0, \pi\}$).

\quad\\
\textbf{Case 1: Different real poles}

In this case we consider $\rho_1 \neq \rho_2$, $\theta_i^0 \in \{0, \pi\}$ and $\beta_i=1$ for $i=1,2$. Then
\begin{eqnarray*}
F(\lambda) & = & \int_{0}^{1} \frac{\sigma^2}{|1-\rho e^{i(\lambda-\theta_1^0)}|^2}dR_1(\rho)
\int_{0}^{1} \frac{\sigma^2}{|1-\rho e^{i(\lambda-\theta^0_2)}|^2}dR_2(\rho) \; .\\
H(\lambda) & = & \int_{0}^{1} \frac{\sigma}{(1-\rho
e^{i(\lambda-\theta_1^0)})}dR_1(\rho)\int_{0}^{1} \frac{\sigma}{(1-\rho
e^{i(\lambda-\theta_2^0)})} dR_2(\rho)\; .
\end{eqnarray*}

Taking, $dR_i(\rho)$ as in \eqref{ec.7.3.1} and applying Lemma \ref{lemma01} we can verify that
\begin{enumerate}

\item [i.] If $-1<d_1, d_2<1$,  and $\theta_1^0 \neq \theta_2^0$, then near $\lambda = \theta_i^0$, for $i=1,2$, $F(\lambda) \cong \frac{1}{|\lambda - \theta_i^0|^{1-d_i}}.$
\item [i'.] If $-1<d_1, d_2<1$,  and $\theta_1^0 = \theta_2^0$, then near $\lambda = \theta_1^0$,
$F(\lambda) \cong \frac{1}{|\lambda - \theta_1^0|^{2-d_1-d_2}}.$
\item [ii.] If $-1<d_1, d_2<0$, and $\theta_1^0 \neq \theta_2^0$, then near $\lambda = \theta_i^0$, for $i=1,2$, $|H(\lambda)| \cong \frac{1}{|\lambda - \theta_i^0|^{-d_i}}.$
\item [ii'.]If $-1<d_1, d_2<0$, and $\theta_1^0 = \theta_2^0$, then near $\lambda = \theta_i^0$, for $i=1,2$, $|H(\lambda)| \cong \frac{1}{|\lambda - \theta_1^0|^{-d_1-d_2}} .$
\end{enumerate}

On the other hand, when $\theta_1^0 \neq \theta_2^0$ we have that $\int F(\lambda) d\lambda < \infty$ if and only if $d_1 >0, d_2 >0$ and $\int |H|^2(\lambda) d\lambda < \infty$ if and only if $d_1 >-\frac{1}{2}, d_2 >-\frac{1}{2}$.

When $\theta_1^0 = \theta_2^0$, $\int F(\lambda) d\lambda < \infty$ if and only if $d_1 + d_2 >1$ and $\int |H|^2(\lambda) d\lambda < \infty$ if and only if $d_1+ d_2 >-1$.

\quad\\
\textbf{Case 2: Complex conjugated poles}

In this case we consider two complex conjugated poles $\rho e^{i\theta}$, $\rho e^{-i\theta} $.
\begin{eqnarray*}
F(\lambda) & = & \int_{0}^{1}\int_{-\pi}^{\pi} \frac{\sigma^2}{|1-\rho e^{i(\lambda-\theta)}|^2
|1-\rho e^{i(\lambda+\theta)}|^2} dR(\rho) dQ(\theta) \; .\\
H(\lambda) & = & \int_{0}^{1}\int_{-\pi}^{\pi} \frac{\sigma^2}{(1-\rho e^{i(\lambda-\theta)})
(1-\rho e^{i(\lambda+\theta)})} dR(\rho) dQ(\theta) \; .
\end{eqnarray*}

We take $dR(\rho)$ as in \eqref{ec.7.3.1} and $Q$ as in \eqref{ec.7.3.2}. Then, applying Lemma~\ref{lemma01} and  Lemma~\ref{lemma02} we obtain

\begin{enumerate}

\item [i.] If $-1<d<2+\beta$, $\beta\leq 1$  and $\theta^0 \in \{0, \pi\}$, then near $\lambda = \theta^0$
$F(\lambda) \cong \frac{1}{|\lambda - \theta_1^0|^{2+ \beta - d}}.$
\item [i'.] If $-1<d<\beta $, $\beta \leq 1$  and $\theta^0 \notin \{0, \pi\}$, then near $\lambda = \theta^0$
$F(\lambda) \cong \frac{1}{|\lambda \mp \theta^0|^{\beta-d}}.$
\item [ii.] If $-1<d<\beta$, $\beta\leq 1$  and $\theta^0 \in \{0, \pi\}$, then near $\lambda = \theta^0$
$|H(\lambda)| \cong \frac{1}{|\lambda - \theta_0|^{\beta -d}}.$
\item [ii'.] If $-1<d<\beta -1 $, $\beta \leq 1$  and $\theta^0 \notin \{0, \pi\}$, then near $\lambda = \pm \theta^0$
$|H|(\lambda)\cong \frac{1}{|\lambda \mp \theta^0|^{-1 + \beta-d}} .$
\end{enumerate}

On the other hand, when $\theta^0 \in \{0,\pi\}$ we have that $\int F(\lambda) d \lambda < \infty$ if and only if $d > 1 + \beta$ and $\int |H|^2(\lambda) d \lambda < \infty$ if and only if $d >\beta - \frac{1}{2}$.

When $\theta^0 \notin \{0,\pi\}$, then $\int F(\lambda) d \lambda < \infty$ if and only if $d>-1+ \beta$ and $\int |H|^2(\lambda) d \lambda < \infty$ if and only if $d >\beta- \frac{3}{2}$.

Finally, we can summarize these results in the following theorem.

\begin{theorem}\label{theo.7.3.2}[ Aggregation of $AR(2)$ processes and long memory.]
If we consider the aggregation of $AR(2)$ processes with random parameters satisfying Assumption~B1, then we have
\begin{enumerate}

\item \textbf{Independent innovations case:}\\
      Different real poles: $\beta_i=1$, $\theta_i \in \{0,\pi\}$ for $i=1,2$.
      \begin{itemize}
      \item $\theta_1^0= \theta_2^0$, $X$ exists if and only if $d_1+d_2>1$, and it is a LM process if and only if $d_1, d_2<1$.
      \item $\theta_1^0 \neq \theta_2^0$, $X$ exists if and only if  $d_1, d_2> 0$, and it is a LM process if and only if  $min\{d_1, d_2\} <1 $.
      \end{itemize}
      Complex conjugated poles: $\beta \leq 1$.
      \begin{itemize}
      \item $\theta^0 \in \{0, \pi\}$, $X$ exists if and only if  $d>1 + \beta$, and it is a LM process if and only if  $d <2+ \beta$.
      \item $\theta^0 \notin \{0, \pi\}$, $X$ exists if and only if  $d> - 1 + \beta$, and it is a LM process if and only if  $d <\beta$.
      \end{itemize}

\item \textbf{Common innovation case:}\\
      Different real poles: $\beta_i=1$, $\theta_i \in \{0,\pi\}$ for $i=1,2$.
      \begin{itemize}
      \item $\theta_1^0= \theta_2^0$, $X$ exists if and only if  $d_1+d_2>-1$, and it is a LM process if and only if $d_1, d_2<0$.
      \item $\theta_1^0 \neq \theta_2^0$, $X$ exists if and only if  $d_1, d_2> -\frac{1}{2}$, and it is a LM process if and only if  $min\{d_1, d_2\} <0 $.
      \end{itemize}
      Complex conjugated poles: $\beta \leq 1$.
      \begin{itemize}
      \item $\theta^0 \in \{0, \pi\}$, $X$ exists if and only if  $d> \beta -\frac{1}{2}$, and it is a LM process if and only if  $d < \beta$.
      \item $\theta^0 \notin \{0, \pi\}$, $X$ exists if and only if  $d>  \beta -\frac{3}{2}$, and it is a LM process if and only if  $d<\beta -1$.
      \end{itemize}
\end{enumerate}
\end{theorem}

\begin{remark}
For different real random poles with distribution concentrate enough near $\rho=1$ we find the same long memory behavior in
$\mathbb{M}(\mathbb{AR}(2))$ as for $\mathbb{M}(\mathbb{AR}(1))$.
\end{remark}


\subsubsection{Opposite phenomena: disappearance of long memory by randomness of $\theta$ parameter \label{sect.7.3.2.1}}

\hspace{0.4cm} Now we illustrate how the long memory can "disappear" by randomness of $\theta$. Here we only consider the case of complex and not real random poles $\rho e^{i\theta}$, $\rho e^{-i\theta}$ and independent innovations. This result can be easily generalized to the case of interactive innovations. We denote by $\bar{S}_{Q}$ the closed support of $Q$.

As we have already mentioned when $\beta=1$ we consider that the measure $Q$ is a Dirac's delta; i.e. $Q_k(\theta)=\delta(\theta-\theta^0)$. For $0<\beta<1$ the measure $dQ(\theta)= \psi(\theta)|\theta-\theta^0|^{-\beta}d\theta$ is strongly concentrated near $\theta^0$ and for $\beta \leq 0$ this measure is regular. In the particulary case of $\beta=0$ we consider that $\psi$ is such that $Q$ is a diffuse measure. For instance, if we take $\psi(\theta)= \1_{(\tau_1, \tau_2)}(\theta)$ then $Q$ is a uniform measure on $(\tau_1, \tau_2)$. With respect to $dR$, we have that when $-1< d < 0$ this measure is strongly concentrated near
$\rho=1$, and when $d\geq 0$ then $dR$ is a regular measure.

\begin{figure}[h!]
\centering \vspace{-0.4cm}\includegraphics[width=11cm]{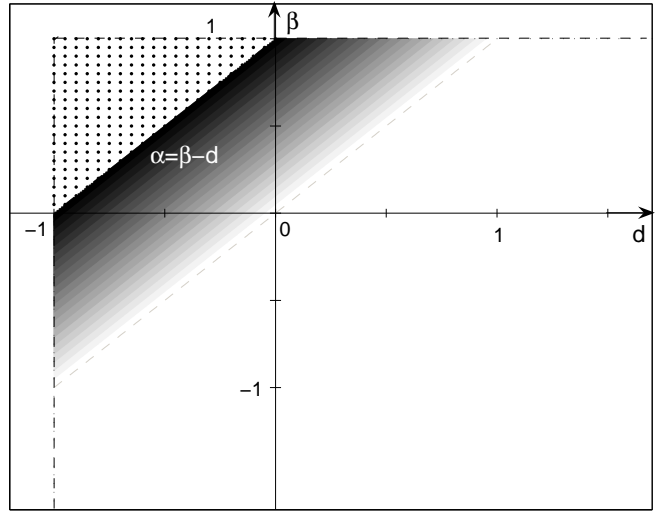}\vspace{-1.2cm} \caption{Values of
$\alpha$-LM parameter (case: $\bar{S}_{Q}$ does not intersect $\{0, \pi\}$). } \label{fig.1}
\end{figure}

\begin{figure}[h!]
\centering \vspace{-0.4cm}\includegraphics[width=11cm]{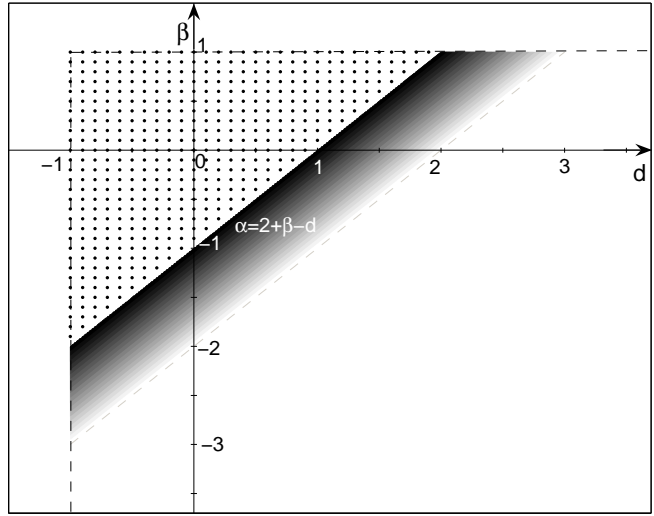} \vspace{-1.2cm} \caption{Values
of $\alpha$-LM parameter (case:  $0$ or  $\pi$ in $\bar{S}_{Q}$).} \label{fig.2}
\end{figure}

In Theorem~\ref{theo.7.3.2} we have given conditions under parameters $\beta$ and $d$ in order to obtain the existence and the long memory property of the aggregation process by means of the aggregation of $AR(2)$ processes. We resume this result in the case of random poles and independent innovations in Figure~\ref{fig.1} and Figure~\ref{fig.2}, where we show the values of $\alpha$-LM
parameter obtained.

Figure \ref{fig.1} represents the existence of the aggregation and the LM property in terms of the values of $(d, \beta)$ in the case where $\bar{S}_{Q}$ does not intersect $\{0, \pi\}$. We recall that we always consider $d\geq-1$ and $\beta\leq 1$. The doted region corresponds to the values of the parameter $(d, \beta)$ for which the aggregation does not exist. The white region corresponds
to values where the aggregation exists but there is no long memory. The remaining region corresponds to the values  for which we obtain the existence and $LM$ property. In this last region we plot the $\alpha$-LM parameter, which is given by $\alpha= \beta-d$. Here, black stands for the maximum value of $\alpha$ (in our case $\alpha=1$, corresponding to an aggregation process strongly dependent) and white stands for its minimum value ($\alpha=0$, coinciding with the case where the aggregation exists but there is no long memory).

We can appreciate that for a value of $d$ fixed, if $\beta$ is small enough, then the $LM$ "disappears", i.e. we can take a measure $Q$ regular enough such that we do not get $LM$ property. For instance, for $0<d<1$ fixed, if we take $Q$ as a Lebesgue measure (i.e. $\beta=0$) such that $\bar{S}_{Q}$ does not intersect $\{0, \pi\}$, then we do not obtain the long memory.

Roughly speaking, it is not enough to have the mixture probability concentrated near $\delta D= \{\rho = 1 \}$ to reach $LM$ by aggregation. But if we take a measure $Q$ such that it is close to a probability with support of Lebesgue measure $0$, for instance $\beta>0$, then we obtain the $LM$. On the other hand, if the measure $dR$ is very concentrated near $\rho=1$, i.e. $d \rightarrow -1$, then it is necessary to take a regular measure $Q$ ($\beta<0$) so that we get the existence of aggregation.

When $\bar{S}_{Q}$ intersects the set $\{0, \pi\}$, we have that $\alpha=\beta - (d-2)$, and we obtain a similar behaviour to the precedent case. We note that in this case the corresponding graphic is the same but translated two units in $d$-axis, see Figure \ref{fig.2}.

\subsection{General case of $AR(p)$ processes \label{sect.7.3.3}}

\hspace{0.4cm} Now, we present the general result in the case of $AR(p)$ processes that allows us to give the condition under the measures $dR$ and $dQ$ in order to obtain  the existence and the long memory of the aggregation process. This result is based on Lemma~\ref{lemma01} and Lemma~\ref{lemma02}.

We denote $n_k = 2m_k$ in the case of independent innovations and $n_k= m_k$ in the common innovation case.

\begin{theorem}\label{theo.7.3.1.1}[ Aggregation of AR(p) processes and long memory.]
Let  $\{Z^i:\,  i\in  \mathbb{N}\}$  be  a sequence  of $AR(p)$ processes  defined  as in \eqref{ec.7.2.0},  let $\{\rho^j \}$ and $\{\theta^j\}$  be the corresponding sequences of random parameter vectors, with $\rho^j=(\rho_1^j,..., \rho_q^j)$ and $\theta^j=(\theta_1^j,...,\theta_q^j)$  for $j  \in \mathbb{N}$. Let $\rho$  and   $\theta$ satisfying  Assumption~B1. Then, the aggregation  exists if and only if the three following conditions hold:
\begin{itemize}
\item $\sum_{k\in K_1} n_k(1+\1_{\{\beta<1\}}(\beta_k))-2+\beta_k-d_k <1.$

\item $\sum_{k\in K_2} n_k(1+\1_{\{\beta<1\}}(\beta_k))-2+\beta_k-d_k <1.$

\item $\min_{k \in K_3} \{n_k - 2 +\beta_k - d_k\}< 1$
\end{itemize}
where $K_1=\{k:  \theta_k^0 =0 \}$, $K_2=\{k:  \theta_k^0 =\pi \}$ and $K_3=\{k: \theta_k^0 \notin
\{0,\pi\} \}$.
\quad\\
Moreover, $X$ is a long memory process if and only if some of following conditions is satisfied:
\begin{itemize}

\item $\forall k\in K_1$,  $n_k \1_{\{\beta<1\}}(\beta_k)-1<d_k<n_k(1+1_{\{\beta<1\}}
(\beta_k))-2+\beta_k$.

\item $\forall k\in K_2$, $n_k \1_{\{\beta<1\}}(\beta_k)-1<d_k<n_k(1+1_{\{\beta<1\}}
(\beta_k))-2+\beta_k$.

\item There exists $ k\in K_3$ such that  $n_k-2<d_k<n_k-2+\beta_k$.
\end{itemize}
\end{theorem}

The   proof   of  this   theorem   is   similar to the proof that we will give for the continuous case.



\section{Aggregation of $OU(p)$ processes and long memory \label{sect.7.4}}

\hspace{0.4cm} The purpose of this sections is to establish the analogues results to Section~\ref{sect.7.2} for the continuous time case.

We consider elementary processes $Z^i$  as stationary solutions of  $p$-order linear stochastic differential equations $LSDE(p)$ driven  by a standard Brownian motion $W^j$. We call these elementary processes Ornstein-Uhlenbeck processes of order $p$, $OU(p)$.

We consider that the characteristic polynomial associated to $LSDE(p)$ equation is factorized as
\begin{equation*}
A(s, \y)=\prod_{k=1}^p (s + \y_k)
\end{equation*}
where, $\y=(\y_1, \ldots, \y_p)$ is the vector of random roots, whose distribution has support $(0,\infty)^p$. Let $c_s(\y)$ be the inverse Laplace transform of $A(s,\y)^{-1}$ defined on $\mathbb{R^{+}}$. We define the $OU(p)$ process  with characteristic polynomial $A(s,\y)$ by
\begin{equation} \label{ec.7.5.3}
Z_t(\y):= \int_{-\infty}^t c_s(\y) dW_s,\quad t\in \mathbb{R}^+.
\end{equation}
This process is, $\mu-a.s.$, a stationary centered Gaussian process with spectral density $g(\lambda,\y)= \sigma^2|A(s,\y)|^{-2}$. We consider the transfer function $h(\lambda,\y)= \sigma A(s,\y)^{-1}$.

If $F$  is the mixture  given  by $F(\lambda)=\E\left[g(\lambda,\y)\right]$, then the process given in \eqref{ec.7.5.3} is well defined, when $\y$ is a random vector, if the distribution $\mu$ of $\y$ satisfies condition~\eqref{ec.7.1.2}.

Let us consider a sequence $\{Z^i_t: \, i \in \N\}$ of $OU(p)$ processes and note by $\{W^i\}$ and $Y=\{\y^i\}$ the corresponding sequences of Gaussian innovations and random parameter vectors. We assume that $\{W^i\}$ and $Y$ satisfy the Assumption~A1 and that interaction $\chi$ satisfies the hypothesis of Theorem 1 given in \cite{Dacunha&Fermin.L.Notes}.

We denote by $\mathbb{OU}(p)$ the class of $OU(p)$ processes and by $\mathbb{M}(\mathbb{OU}(p))$ the class of processes that can be obtained by aggregation of elementary processes in $\mathbb{OU}(p)$.


In the sequel, we consider that
$$A(s,\y)=\prod_{k=1}^n(s+r_k)^{m_k}\prod_{k=n+1}^{q}[(s+r_k+i\tau_k) (s+r_k-i\tau_k)]^{m_k},$$
is a polynomial  of grade  $p$ with $n$  real roots and  $2(q-n)$ complex pairwise conjugate  roots having  strictly positive real parts, where $n\leq q$,  $m_k$ is the  multiplicity of the roots $\y_k = r_k\pm i\tau_k$, with      $\tau_k=0$      for       $1\leq k\leq n$,  and $p=\sum_{k=1}^{n}m_k+2\sum_{j=n+1}^{q}m_j$.

For $r=(r_1,...,r_q)$ and $\tau=(\tau_1,...,\tau_q)$ fixed, under condition~\eqref{ec.7.1.2} $Z_t$ is a stationary centered Gaussian process with spectral density
\begin{equation}\label{ec.7.6.1}
g(\lambda,r,\tau)=\prod_{k=1}^{n}\frac{1}{(\lambda^2+r_k^2)^{m_k}} \prod_{k=n+1}^{q}
\frac{1}{[\left((\lambda-\tau_k)^2+r_k^2\right) \left((\lambda+\tau_k)^2+r_k^2\right)]^{m_k}},
\end{equation}
and transfer function
\begin{equation}\label{ec.7.6.2}
h(\lambda,r,\tau)=\prod_{k=1}^{n}\frac{1}{(i \lambda + r_k)^{m_k}} \prod_{k=n+1}^{q}
\frac{1}{[\left((\lambda-\tau_k)i + r_k\right) \left((\lambda+\tau_k)i + r_k\right)]^{m_k}}.
\end{equation}

As in the discrete case, we will consider the following assumption.

\quad\\
\textbf{Assumption B2:} \textit{Let $r_1,..., r_q$ be independent  random variables, such that
$r_k$ has distribution
\begin{equation}\label{ec.7.6.3}
dR_k(r)=|r|^{d_k} \varphi_k(r)dr,
\end{equation}
where  $\varphi_k$ is a  bounded positive  function, continuous  in $r=0$ with   $\varphi_k(0)>0$.
Let   $\tau_1,...,\tau_q$ be independent random variables and  independent of $r_1,...,  r_q$, such
that $\tau_k$ has distribution
\begin{equation}\label{ec.7.6.4}
dQ_k(\tau)=\frac{\psi_k(\tau)d\tau}{|\tau-\tau_k^0|^{\beta_k}},
\end{equation}
where  $\beta_k  \leq  1$  $\psi_k$  is  a  bounded  positive function, continuous  in
$\tau=\tau_k^0$  and  $\psi_k(\tau_k^0)>0$. By convention when $\beta_k=1$ we will consider
$Q_k(\tau)=\delta(\tau-\tau_k^0)$ and $\tau_k^0 = 0$ for $1\leq k\leq n$.}
\quad\\

Let us now study the local behavior of the spectral densities $F$ and $|H|^2$, defined by equations \eqref{ec.7.1.1} and \eqref{ec.7.1.3} respectively. We take $n_k=2m_k$ in the case of independent innovations and $n_k= m_k$ in the common innovation case. Then, under Assumption~B2 we obtain the following technical lemmas.

\begin{lemma}\label{lemma1} \quad
\begin{enumerate}

\item If $-1<d_k<n_k-1$, then near $\lambda=0$
\begin{eqnarray*}
\int_0^{\infty} \frac{|r|^{d_k}\varphi_k(r)}{(\lambda^2+r^2)^{n_k}/2}dr \sim
\frac{\varphi_k(0)}{|\lambda|^{n_k-1-d_k}}\int_0^{\infty}\frac{u^{d_k}du} {(1+u^2)^{n_k/2}}.
\end{eqnarray*}

\item If $-1<d_k<n_k-1$ and $\tau_k^0 \neq 0$, then near $\lambda =\pm
\tau_k^0$
\begin{eqnarray*}
\int_0^{\infty}\!\!\frac{|r|^{d_k}\varphi_k(r)dr}
{[((\lambda-\tau_k^0)^2+r^2)((\lambda+\tau_k^0)^2+r^2)]^{n_k/2}} \sim
\frac{\varphi_k(0)(2\tau_k^0)^{-n_k}}{|\lambda\mp\tau_k^0|^{n_k-1-d_k}}
\!\int_0^{\infty}\!\!\!\frac{u^{d_k} du}{(1+u^2)^{n_k/2}}.
\end{eqnarray*}
\end{enumerate}
\end{lemma}
\vspace{0.4cm}

\begin{proof}\qquad
The first part  is shown by making the  variable change $r=|\lambda|u$ and  then  taking limits  as
$\lambda\rightarrow 0$. The  condition $-1<d_k<n_k-1$ implies the  convergence  of the following
integral $\int_0^{\infty}\frac{u^{d_k}}{(1+u^2)^{n_k/2}}du$.

The second part can be proved in the same way, with the variable change $r=|\lambda \mp \tau_k^0|u$
and  taking limits as $\lambda \rightarrow \pm \tau_k^0$.

\end{proof}

\vspace{0.4cm}

\begin{lemma}\label{lemma2} Let $\alpha_k<1$ and
$$f_k(\lambda) = \int_{-\infty}^{\infty}\int_0^{\infty}
\frac{|r|^{d_k}\varphi_k(r)\psi_k(\tau)|\tau-\tau_k^0|^{-\alpha_k}dr d\tau}
{[((\lambda-\tau)^2+r^2)((\lambda+\tau)^2+r^2)]^{n_k/2}}.$$
\begin{enumerate}
\item If $ n_k-1<d_k<2n_k-2+\alpha_k$ and $\tau_k^0 = 0$, then near
$\lambda = 0$
\begin{eqnarray*}
f_k(\lambda) \sim \frac{\varphi_k(0) \psi_k(0)}{|\lambda|^{2n_k-2-d_k+\alpha_k}}
\int_{-\infty}^{\infty}\int_0^{\infty} \frac{u^{d_k} |\theta|^{-\alpha_k}du
d\theta}{[((\theta-1)^2+u^2)((\theta+1)^2+u^2)]^{n_k/2}}.
\end{eqnarray*}
\item If $n_k-2<d_k<n_k-2+\alpha_k$ and $\tau_k^0 \neq 0$, then near
$\lambda = \pm \tau_k^0$
\begin{eqnarray*}
f_k(\lambda) \sim \frac{\varphi_k(0) \psi_k(\tau_k^0)(2\tau_k^0)^{-n_k}}{|\lambda \mp
\tau_k^0|^{n_k-2+\alpha_k-d_k}}\!
\int_{-\infty}^{\infty}\!\frac{|\theta|^{-\alpha_k}d\theta}{|\theta \mp 1|^{n_k-1-d_k}}
\!\int_0^{\infty} \!  \frac{u^{d_k} du}{ (1+u^2)^{n_k/2}}.
\end{eqnarray*}
\end{enumerate}
\end{lemma}

\vspace{0.4cm}

\begin{proof}\qquad
We show point $1$ by making the variable changes $r=|\lambda|u$ and $\tau=\lambda\theta$ and then taking limits as $\lambda \rightarrow 0$. The result holds if the integral
$$I:=\int_{-\infty}^{\infty} \int_0^{\infty}\frac{u^{d_k} du
d\theta}{[((\theta-1)^2+u^2)((\theta+1)^2+u^2)]^{n_k/2}|\theta|^{\alpha_k} }$$
is convergent. To verify that we will use that
$$u^2(1+\theta^2+u^2) <((\theta-1)^2+u^2)((\theta+1)^2+u^2)$$
and   then    we   will   make   the    following   variable change: $u=(1+\theta^2)^{\frac{1}{2}}r$, from where
$$I \leq \int_{-\infty}^{\infty}\int_0^{\infty}\frac{u^{d_k-n_k}du
d\theta}{ (1+ \theta^2 + u^2)^{n_k/2}|\theta|^{\alpha_k} } =
\int_0^{\infty}\frac{2\theta^{-\alpha_k}d\theta}{(1+\theta^2)^{\frac{2n_k-1-d_k}{2}}}
\int_0^{\infty} \frac{r^{d_k-n_k} dr}{(1+r^2)^{n_k/2}}.$$
Finally we have that these two last integrals converge if $\alpha_k<1$ and $n_k-1<d_k<2n_k-2+\alpha_k$.

To prove point $2$ we make two variable changes $r=|\lambda \mp \tau|u$ and $\tau-\tau_k^0=(\lambda \mp \tau_k^0)\theta$, which give us
\begin{eqnarray*}
f_k(\lambda) & = & \frac{1}{|\lambda \mp \tau_k^0|^{(n_k-2)+\alpha_k-d_k}}\int_{-\infty}^{\infty}
\frac{\psi_k(\tau_k^0+(\lambda \mp \tau_k^0)\theta)} {|\theta \mp
1|^{n_k-1-d_k}|\theta|^{\alpha_k}}\\ & & \times \int_0^{\infty}\frac{u^{d_k}\varphi_k(|\lambda \mp
\tau_k^0||\theta \mp 1|u) du d\theta} {[(1+u^2)([(\lambda \pm \tau_k^0) \pm (\lambda \mp
\tau_k^0)\theta]^2+[(\lambda \mp \tau_k^0)(\theta \mp 1)u]^2)]^{n_k/2}}.
\end{eqnarray*}
Then, taking limits as $\lambda \rightarrow \pm \tau_k^0$,
$$f_k(\lambda) \sim \frac{\varphi_k(0)\psi_k(\tau_k^0)(2{\tau_k^0})^{-n_k}
}{|\lambda \mp \tau_k^0|^{n_k-2+\alpha_k-d_k}} \int_{-\infty}^{\infty} \frac{d\theta}{|\theta
\mp1|^{n_k-1-d_k}|\theta|^{\alpha_k}} \int_0^{\infty}\frac{u^{d_k}du}{(1+u^2)^{n_k/2}}.$$
Finally, we can see that  these two integrals converge if $\alpha_k<1$ and $n_k-2<d_k<n_k-2+\alpha_k$.

\end{proof}

\vspace{0.4cm}


From Lemma~\ref{lemma1} and Lemma~\ref{lemma2} and under Assumption~B2 we show the following theorem which gives the condition over parameters $d, \beta$ that allows us to obtain the long memory for the aggregation of $OU(p)$ processes.

\vspace{0.4cm}

\begin{theorem}\label{theo.7.5.1}[ Aggregation of $OU(p)$ processes and long memory.]
Let  $\{Z^i:\,  i\in  \mathbb{N}\}$  be  an i.i.d  sequence  of $OU(p)$ processes  defined as  in \eqref{ec.7.5.3},  let $\{W^i\}$, $\{r^i\}$ and $\{\tau^i\}$ be the corresponding  i.i.d sequences of Brownian motions and of  random parameters vectors.  Let $r$ and $\tau$ be random vector satisfying Assumption~B2. Then,  the aggregation  exists if condition \eqref{ec.7.1.2} holds  and
there is long memory if and only if $$d_k<n_k(1+\1_{\{0\}}(\tau_k^0) \1_{\{\beta<1\}}(\beta_k))-2+\beta_k, \quad  for  \,   some  \quad k\in\{1,...,q\}.$$
\end{theorem}

\vspace{0.4cm}

\begin{proof}\qquad It is clear that condition \eqref{ec.7.1.2} implies the aggregation existence. On the  other hand, as long memory of the aggregation is  related to the spectral density singularities, then by Lemma~\ref{lemma1} and Lemma~\ref{lemma2}, the result
holds.

\end{proof}

\vspace{0.4cm}

This theorem generalizes the results given by other authors in the following sense: in our approach $\tau$ is a  random parameters  vector  and  not a constant vector, moreover we consider multiple roots. In the  case where $\beta_k=1$ for $1\leq k \leq q$, i.e. when we  consider the parameters vector $\tau$ fixed, the necessary and sufficient conditions for the aggregation existence can be easily  written in terms of the parameters  $d_k$ and $n_k$, as we can see  in the Corollary \ref{cor.7.6.1}. In  the general case, the aggregation existence does not  depend  on  the decay of functions $\varphi_k$ and $\psi_k$, for $1\leq k \leq q$. For instance, it can be seen that in the case of $OU(2)$ processes if we consider $\varphi_k$ bounded and $\psi_k(\tau)\sim|\tau|^{-\beta_k}\wedge|\tau|^{-\beta_k^{'}}$ with $\beta_k <1$  and $\beta_k^{'}>0$,  then  the  aggregation exists if and only if  $1+\beta_k < d_k < 2$.

The following  corollary generalize the  result given in \cite{Oppenheim&Viano2004}, because these authors consider only the case of simple roots and $\tau$ fixed.

\begin{corollary}\label{cor.7.6.1}
Let  $\{Z^i: \,  i\in  \mathbb{N}\}$  be  an i.i.d  sequence  of $OU(p)$ processes, $\{W^i\}$, $\{r^i\}$ and $\{\tau^i\}$ the corresponding i.i.d sequences of Brownian motions and of parameters vectors. If  we consider $\alpha_k=1$ for $1\leq k \leq q$ and the parameters  vectors $r$ and $\tau$ satisfying Assumption~B2, then the  aggregation exists  and there is long  memory if and
only if  the following conditions hold:
\begin{itemize}
\item $-1<d_k<n_k-1$ for $1\leq k \leq n$.

\item $n_k-2<d_k<n_k-1$, for $n< k \leq q$.

\item $\sum_{k=1}^n n_k-d_k<n+1$.
\end{itemize}
\end{corollary}
\vspace{0.4cm}
\begin{proof}\qquad By applying Lemma~\ref{lemma1} and  if $-1<d_k<n_k-1$, for $1\leq k \leq q$, then near to $\lambda=0$ we have that the spectral density $G$ of the aggregation process is such that
$$G(\lambda)\sim\prod_{k=1}^n\frac{1}{|\lambda|^{n_k-1-d_k}}\varphi_k(0) \int_0^{\infty}\frac{u^{d_k}}{(1+u^2)^{\frac{n_k}{2}}}du,$$
and near to $\lambda =\pm \tau_k^0$
$$G(\lambda)\sim\frac{1}{|\lambda\mp\tau_k^0|^{n_k-1-d_k}} \frac{\varphi_k(0)}{(2\tau_k^0)^{n_k}}\int_0^{\infty}\frac{u^{d_k} du} {(1+u^2)^{\frac{n_k}{2}}}.$$
Moreover, we have that
$$G(\lambda) \leq \frac{1}{\prod_{k=1}^n \lambda^{n_k}\prod_{k=n+1}^q(\lambda^2-{\tau_k^0}^2)^{n_k}},$$
which allows us  to bound $G$ when $\lambda \rightarrow \pm \infty$. Then,  the aggregation  exists if and only if  $-1<d_k<n_k-1$ for $1\leq k \leq q$, $n_k-1-d_k<1$ for  $n < k \leq q$, and $\sum_{k=1}^n n_k-1-d_k<1$.  Furthermore, there is long memory  if and only if  $n_k-1-d_k>0$ for some $k \in\{1,...,q\}$, from where the theorem holds.

\end{proof}
\vspace{0.4cm}

A slight modification in the proof of Theorem~\ref{theo.7.5.1} allows us to extend these results when
we consider the measures $dQ_k$ as
$$dQ_k(\tau)=\sum_{j=1}^{n_k}p_j\delta(\tau-\tau_k^j)+\Psi_k(\tau)d\tau$$
where $\Psi_k$  is a positive function with  singular points $s_1,\ldots, s_{l_k}$ such that  for each $s_j$ there exists  a function $\psi_{j,k}$, regular in a neighborhood $V(s_j)$ of $s_j$, such that
$$\Psi_k(s)\sim  \frac{\psi_{j,k}(s)}{|s-s_j|^{\beta_{j,k}}}, \quad  \mbox{ for }  \quad s \in V(s_j)$$
and $\Psi_k$ a function bounded out of $\bigcup_j V(s_j)$.

\begin{remark}
In the  aggregation of $OU$ processes the phenomenon of disappearance of long memory also can
happen  by randomness of parameter $\tau$. The analysis is similar to the one given for the
discrete time case.
\end{remark}


\section*{Acknowledgements} We thank to E. Moulines (TELECOM Paris-Tech) for helpful comments.


\bibliographystyle{acm}
\bibliography{biblio_aggreg_LM}

\begin{thebibliography}{10}

\bibitem{Barndorff}
{\sc Barndorff-Nielsen, O.~E.}
\newblock {Superposition of Ornstein-Uhlenbeck type processes}.
\newblock {\em { Theory of Probability and its Applications} 45\/} (2001),
  175--194.

\bibitem{Beran}
{\sc Beran, J.}
\newblock {\em {Statistics for long-memory processes }}.
\newblock Chapman \& Hall, New York, (1994).

\bibitem{Celov}
{\sc Celov, D., Leipus, R., and Philippe, A.}
\newblock {Time series aggregation, disaggregation and long memory}.
\newblock {\em {Preprint, arXiv:math/0702821v1 [math.ST]}\/}.

\bibitem{Dacunha&Fermin.L.Notes}
{\sc Dacunha-Castelle, D., and Ferm\'{\i}n, L.}
\newblock {Aggregation of Doubly Stochastic Interactive Gaussian Processes and
  Toeplitz forms of $U$-Statistics}.
\newblock {\em {In Dependence in Probability and Statistics, Series: Lecture
  Notes in Statistic.} 187\/} (2006).

\bibitem{Fermin}
{\sc Ferm\'{\i}n, L.}
\newblock {Aggregation of weakly dependent doubly stochastic processes}.
\newblock {\em {arXiv:0805.1949v1 [math.PR].}\/} (2008).

\bibitem{Goncalvez1988}
{\sc Gon\c{c}alvez, E., and Gourieroux, C.}
\newblock {Agr\'egation de processus autorégressifs d'ordre 1}.
\newblock {\em { Annales d'Economie et de Statistique } 12\/} (1988), 127--149.

\bibitem{Granger}
{\sc Granger, C.}
\newblock {Long Memory relationships and the aggregate of dinamic models}.
\newblock {\em {Journal of Econometrics} 14\/} (1980), 227--238.

\bibitem{Igloi&Terdik}
{\sc Igloi, E., and Terdik, G.}
\newblock {Long-range dependence through gamma-mixed Ornstein-Uhlenbeck process
  }.
\newblock {\em {Electronic Journal of Probability} 4\/} (1999), 1--33.

\bibitem{Linden}
{\sc Linden, M.}
\newblock {Time series properties of aggregated AR(1) processes with uniformly
  distributed coefficients }.
\newblock {\em {Economics Letters} 64\/} (1999), 31--36.

\bibitem{Lippi&Zaffaroni}
{\sc Lippi, M., and Zaffaroni, P.}
\newblock {Aggegation of simple linear dynamics: exact asymptotic results}.
\newblock {\em {Econometrics Discussion Paper 350, STICERD-LSE }\/}.

\bibitem{Nikolai}
{\sc Nikolai, L., and Emanuele, T.}
\newblock {Convergence of integrated superpositions of Ornstein-Uhlenbeck
  processes to fractional Brownian motion.}
\newblock {\em Stochastics 77}, 6 (2005), 477--499.

\bibitem{Oppenheim&Viano2004}
{\sc Oppenheim, G., and Viano, M.-C.}
\newblock {Aggregation of ramdom parameters Ornstein-Uhlenbeck or AR processes:
  some convergence results}.
\newblock {\em {Journal of Time Series Analysis } 25}, {3} (2004), 335--350.

\bibitem{Terence1}
{\sc Terence, T., and to~W., K.}
\newblock {Time series properties of aggregated AR(2) processes}.
\newblock {\em {Economics Letters} 73\/} (2001), 325--332.

\end{thebibliography}


\end{document}